\documentclass[a4paper]{article}
\usepackage{upgreek}
\usepackage{graphicx}
\usepackage{color}
\usepackage{enumitem}
\usepackage{amsmath, amsthm, amsfonts, amssymb, amscd}
\usepackage{setspace}
\usepackage{latexsym}
\usepackage[colorlinks=true,hypertexnames=false,linkcolor=blue,citecolor=blue, backref=page]{hyperref}
\usepackage{mathrsfs}
\usepackage{newclude}
\usepackage{verbatim}
 \usepackage[usenames,dvipsnames]{xcolor}
\usepackage{cite}
\usepackage[noblocks]{authblk}
\usepackage[utf8]{inputenc}
\DeclareMathAlphabet\mathbfcal{OMS}{cmsy}{b}{n}

\advance \textheight by \topskip
\numberwithin{equation}{section}
\usepackage{fancyhdr}
\newenvironment{lista}%
  {\begin{itemize}[label=--,topsep=-8pt, partopsep=0pt]%
    \setlength{\itemsep}{0pt}%
    \setlength{\parskip}{0pt}}%
  {\end{itemize}}%%%%%%%%%%%%%%%%%%%%%%%%%%%%%%%%%%%%%%%%
\newtheorem{defin}{Definition}[section]
\newtheorem{thm}{Theorem}
\newtheorem{lem}{Lemma}[section]
\newtheorem{prop}[lem]{Proposition}

\makeatletter
\newenvironment{prova}[1][\proofname]{%\par
\pushQED{\qed}%
\normalfont \topsep0\p@\@plus0\p@\relax
\trivlist
\item\relax
{\bfseries
#1\@addpunct{.}}\hspace\labelsep\ignorespaces
}{%
\popQED\endtrivlist\@endpefalse
}
\makeatother
\DeclareMathOperator{\gdc}{gdc}
\def\ackname{Acknowledgements}%
\newenvironment{ack}[1][\ackname]%
{\ifx#1\empty\else\subsection*{#1.}\fi\par}
{\par}
\usepackage{cancel}
\usepackage[normalem]{ulem}
\newcommand\redst{\bgroup\markoverwith{\textcolor{red}{\rule[0.5ex]{2pt}{0.7pt}}}\ULon}
\newcommand\bluest{\bgroup\markoverwith{\textcolor{blue}{\rule[0.5ex]{2pt}{0.7pt}}}\ULon}
\usepackage[textwidth=2.7cm, color=red!20]{todonotes}
\newcounter{nota}[section]

\newcommand{\notaIL}[3][]{\refstepcounter{nota}%
{%
\todo[linecolor=blue!40,backgroundcolor={blue!20!},size=\footnotesize,inline]
{\textbf{{#1}{#2}[\thesection.\thenota]:}~{#3}}%
}}

\title{Symmetric periodic orbits in symmetric billiards}
\author{
Geraldo C\'esar Ferreira Gon\c calves, Sylvie Oliffson Kamphorst and S\^onia Pinto-de-Carvalho\\
}
\date{}
\begin{document}
\maketitle

\begin{abstract}
In this text we study  billiards on ovals and investigate some consequences of a rotational symmetry of the boundary on the dynamics. As it simplifies some calculations, the symmetry helps to obtain the results.  
We focus on periodic orbits with the same symmetry of the boundary which always exist and prove that in general half of them are elliptic and Moser stable and the other half are hyperbolic with homo(hetero)clinic intersections.
\end{abstract}

%\ams{37D30, 37C70}

 %intro
\section{Introduction and general facts about billiards} \label{sec:intro}

%\notaIL{Syok} {para ideia do estilo, olhar os artigos recentes do Bunimovich \cite{bunim22} e Jin-Zhang \cite{jin22} no Nonlinearity 2022 }

The billiard problem
\cite{birk,markarian,tab} is the study of the free motion of a particle a plane region bounded by a curve $\Gamma$.
The particle travels at a constant speed inside the region and undergoes elastic reflections at the impacts with the boundary,
following a polygonal trajectory.
The billiard dynamics is completely described by given at each collision the impact point on the boundary
and the direction of motion immediately after it and thus is given by a two dimensional conservative map. 
The billiard map associates to each  pair (impact point, direction of motion) the next one. 

An oval is a smooth, regular, closed, plane, 
curve with strictly positive curvature $K$. 
Billiards on ovals are known as Birkoff's Billiards. 
For them, it is useful to use  
the angle $\phi$ of the tangent vector with a fixed axis to parametrize the curve.  So the impact point is identified by $\phi$.
The direction of motion after the impact is specified  by the tangential momentum $p=\cos \alpha$, where $\alpha$ is the outgoing angle measured from the tangent vector. 
Using the variables $(\phi,p)$ the  billiard motion in the region bounded by the curve $\Gamma (\phi)$ is then described by a map $T$
of the cylinder
$ \Omega = \mathbb{R}/{2\pi}\mathbb{Z} \times (-1,1)$
which will be represented in our figures as a rectangle. 
To each orbit $\mathcal{O}(\phi_0,p_0) = \{  T^n(\phi_0, p_0) \}$ in the phase space $\Omega$ corresponds an oriented polygonal trajectory in the plane with vertices at the points $\{ \Gamma(\phi_n) \}$ of the boundary.
The billiard map $T$ has some well known properties \cite{markarian}: it is invertible (by reversing the orientation of a trajectory, one obtain the inverse orbit),  is a twist diffeomorphism 
and preserves the measure $ R \, d\phi \, dp$, where $R=1/K$ is the curvature radius.
The inversibility implies that the phase space of the billiard map is symmetric under the reflection
along the horizontal middle line $p=0$.
Moreover, Birkoff's Min-Max implies the existence of many periodic orbits.

In this work, we focus on  billiards in symmetric ovals, as described in Section~\ref{sec:curvas}, and investigate how the symmetry influences the dynamics. We assume that our ovals are $C^\infty$ although most of our proofs also hold for less differentiable curves. 
Our main result is that a "typical" symmetric billiard has both elliptic islands 
and hyperbolic orbits with transverse intersection of the stable and unstable curves. In this sense, symmetric billiards have mixed dynamics.  This is a relevant result, as Birkoff's Min-Max Theorem does not guarantee the existence of elliptic orbits.

The dynamics obtained is associated to orbits with the same symmetry of the curve.  
We construct in  Section~\ref{sec:poligonos} regular polygons inscribed in a $n$-symmetric oval and  show that the polygons of maximum and minimum perimeter are effectively trajectories of the billiard motion. 
The results we obtain rely on the fact that periodic trajectories are associated to the critical points of the support function. So, by identifying symmetric ovals by their support functions, we shift our study to periodic functions.

Section~\ref{sec:linear} contains the analysis of the linear stability of the orbits associated to the regular polygons. 
This analysis is simplified by the equivalence relation defined by the symmetry since we can use the billiard map in a quotient space. We show that the elliptic orbits are located at the (non degenerate) minima of the support function and the hyperbolic at the maxima.
So, in a suitable topology, for a ''typical" symmetric oval, half of the symmetric orbits are elliptic and the other half hyperbolic. 

Analyzing the symmetric elliptic orbits, through the minima of the support function, we are able to prove in Section~\ref{sec:elip} that "typical" symmetric billiards have elliptic islands. To do so,  we obtain an explicit expression for the twist coefficient of a quotient map. The perturbation arguments used in \cite{ilhas} are carried over to the
support function.

In Section~\ref{sec:hyp} we turn our attention to the symmetric hyberbolic orbits. 
We show that a perturbation of the support function,  breaks a tangency to create a transverse intersection of the stable and unstable manifolds as in \cite{periodicas}. 

{Finally  in Section~\ref{sec:final}
%\nota{Syok}{vamos querer fazer conjecturas?}
 we present some figures of the phase space of symmetric billiards which, in addition to display the features predicted by our results,  raise some interesting questions and point to further investigation about the existence of invariant rotational curves beyond those close to the boundary predicted by Lazutkin's Theorem \cite{laz}.}

At last, but not least, we point out that the symmetry is crucial here, as it allows the explicit construction of also symmetric periodic orbits
where the computation of the Birkoff coefficient it drastically simplified. So we can prove the existence islands of higher period.
Even though the calculation is more complicated that the one of period two orbits \cite{ss}, it looks more direct that the one in \cite{bunin-grigo}.
We also have no need to study the second order twist coefficient \cite{jin22}.  

 %intro
\section{Symmetric ovals}\label{sec:curvas}

We will study the billiard on ovals (regular, closed, $C^{\infty}$ plane curves with stricly positive curvature) with rotational symmetry. 
Given an integer $n \ge 2$, we say that an oval $\Gamma$ is $n$-symmetric if it is invariant under a 
$\frac{2\pi}{n}$ rotation around a point ${o}$ called the symmetry center which we can assume is at the origin.
The set of $C^{\infty}$ $n$-symmetric ovals
with symmetry center at the origin and parametrized by the angle $\phi $
between the tangent vector $\sigma = (\cos \phi, \sin \phi)$ and the positive horizontal axis $x$ is denoted by  $\mathcal{N}$.

\begin{figure}[h]
\begin{center}
\includegraphics[width=0.3\hsize]{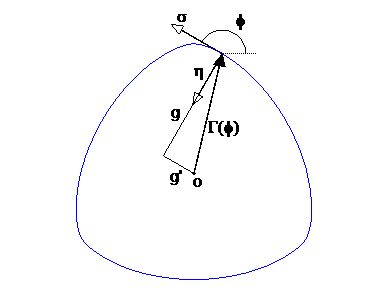}
\end{center}
\caption{Oval $\Gamma$ defined by its support function as in \ref{eq:curva em g}}
\label{fig:coordg}
\end{figure}

The support function of  an oval $\Gamma$ is defined by $g(\phi) = - < \Gamma(\phi), \eta(\phi) >$, where $\eta(\phi) = (-\sin \phi, \cos \phi)$ is the inward normal. 
The positive function $g$ is related to the radius of curvature by $g(\phi) + g''(\phi) = R(\phi) $.
Given a strictly positive, $\frac{2\pi}{n}$ periodic and smooth function $g(\phi)$, such that $g+g'' > 0$,  an $n$-symmetric oval is defined (Figure \ref{fig:coordg}) by
\begin{equation}
\Gamma (\phi)  = - g(\phi) (-\sin \phi, \cos \phi)+ g'(\phi)  (\cos \phi, -\sin \phi)  \label{eq:curva em g}
\end{equation}

Let $S_n= \left. \mathbb{R} \middle / \frac{2\pi }{n}\mathbb{Z} \right.$ 
and consider the set 
$C^{\infty}\left( S_n,\mathbb{R} \right) $ 
with the $C^2$ norm
$\|  . \|_{2}$.
We define
\begin{equation}
\mathcal{G}=\{ g\in C^{\infty}(S_n, \mathbb{R}) ,  g( \phi) >0, g( \phi) +g^{\prime \prime }( \phi )>0 \}
\end{equation}
and refer to it as the set of ({$C^{\infty}$}) support functions. 
 It is clear that $\mathcal{G}$ open in $C^{\infty}(S_n, \mathbb{R})$.

It is easy to verify that we  have a one to one correspondence between the set $\mathcal{N}$ and the set 
$\mathcal{G}$. Using  the topology induced by $\mathcal{G}$ in $\mathcal{N}$
we have then an homeomorphism between these two sets.
In this topology, the set  $\mathcal{N}$ is a metric space with norm
$\left\vert \left\vert \Gamma \right\vert \right\vert_{2}
=\left\vert \left\vert g\right\vert \right\vert_{2}$, where  $g$  is the function in
$\mathcal{G}$ which corresponds to the oval  $\Gamma $ in $\mathcal{N}$.

At this point, it is worthwhile to notice that a translation of the origin does not affect the support function. 
On the other hand, a translation in the argument of the support function, $g(\phi) \to g(\phi - \phi_0)$ results in a rotation of the curve. So, in fact we can associate to a support function a class of (equivalent) curves obtained from the original one by translations and rotations. This clearly does not change the dynamics. So from now on, a curve will be in fact, in our context, an equivalence class. 
Moreover, multiplying the support function by a constant results in a homothety transformation of the curve which also does not affect the dynamics, since our results are based on properties of the critical points of the support function which are scaling invariant.

Finally it will be useful to consider support functions which are Morse functions.
We denote by  $\tilde{\mathcal{G}} \subset {\mathcal{G}}$ 
the subset of $C^{ \infty}\left( S_n,\mathbb{R}\right) $ positive Morse functions $g$ such that $g+g^{\prime\prime} >0$ and refer to it as the set of Morse support functions.
The set $\tilde{\mathcal{N}} \subset \mathcal{N}$ of $C^{\infty}$ $n$-symmetric ovals 
associated to  the set $\tilde{\mathcal{G}}$ is called the set of Morse $n$-symmetric ovals, it is open and dense in $\mathcal{N}$ as is 
$\tilde {\mathcal G}$ in  $\mathcal{G}$.

%%%%%%%%%%%%%%%%%%%%%%%%%%%%%%%%%%%%%%%%%%%%%%%%%%%%%%%

 %curvas
\section{Periodic motion} \label{sec:poligonos}

%%%%%%%%%%%%%%%%%%%%%%%%%%%%%%%%%%%%%%%%%%%%%%%%
The billiard motion in the region bounded by an oval is a polygonal with reflexions on the impacts with the border.
So periodic motion takes place along polygons inscribed in the billiard's boundary and it is clear that an inscribed polygon is a billiard trajectory if and only if its sides make equal angles with the tangent line at the vertices. 
Any oval two has at least two trajectories of period 2, located at the diameters.  More generally, Birkoff's Min-Max theorem \cite{birk} assures the existence of trajectories of any period. Our strategy to find periodic trajectories is to repeat this construction using the symmetry.

Given a $n$-symmetric oval $\Gamma \in   \mathcal{N}$
we construct regular polygons inscribed in it.
\begin{defin}
\label{def:poligonos}
For each given integer $m$ with  $1 \le m \le  n/2  $ 
%\footnote{$\lfloor x\rfloor = floor(x) $  is the largest integer not greater than $x$}
and each $\phi$,  the points
$$
\Gamma(\phi), \Gamma(\phi + \frac{2m\pi}{n}), \Gamma(\phi +  {\frac{4m\pi}{n}}),
 \ldots,
\Gamma(\phi + \frac{2(\tilde{n}-1)m\pi}{n})
$$
where $\tilde n = \frac{n}{\gdc(n,m)}$,
are the vertex of a regular $\tilde n$-polygon denoted by $P_{m}(\phi)$.
\end{defin}

\begin{figure}[h]
\begin{center}
\includegraphics[width=0.2\hsize]{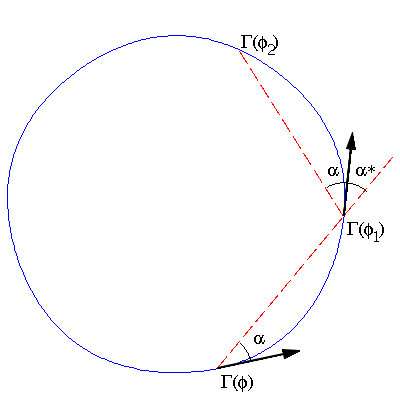}
\hskip 1.8cm
\includegraphics[width=0.25\hsize]{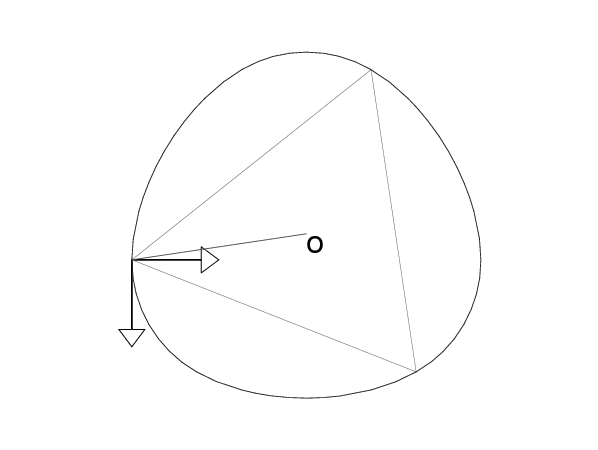}
\includegraphics[width=0.25\hsize]{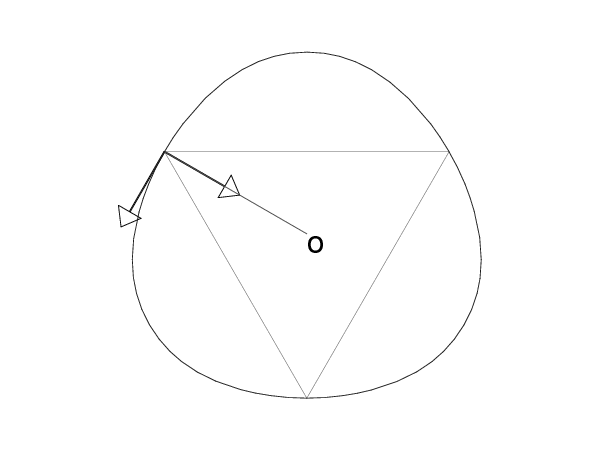}
\caption{Reflexion property (left). 
Equilateral triangles inscribed in a 3-symmetric oval: the one on the center is not a billiard trajectory,
while the one on the right is.}
\label{fig:reflex}
\end{center}
\end{figure}

Clearly $P_{m}(\phi)$ is equal to $P_{m}(\phi+\frac{2m\pi}{n})$ and is congruent to $P_{m}(\phi+\frac{2\pi}{n})$ 
so we have $\gdc(n,m)$ distinct $\tilde n$-polygons.
The length of the sides of a polygon $P_{m}(\phi)$ is given by 
$L_m(\phi) =  ||\Gamma(\phi_1) - \Gamma(\phi)||$ 
where $\phi_k = \phi+\frac{2km\pi}{n}$. It is easy to verify that  
\begin{eqnarray*}
L_m(\phi) \frac{d L_m}{d\phi} &=&
 R(\phi) \left(
 \left< \Gamma\left(\phi_1\right) 
- \Gamma(\phi),\sigma\left(\phi_1\right) \right>
-
\left< \Gamma\left(\phi_2\right) 
- \Gamma\left(\phi_1\right) ,\sigma\left(\phi_1\right) \right>
\right)
\\
&=&
 R(\phi) (\cos \alpha^* - \cos \alpha)
\end{eqnarray*}
an so the angles of the sides at the vertex are equal (the incidence angle $\alpha*$ equals the outgoing angle $\alpha$, as on Figure~\ref{fig:reflex}) 
if and only
if $\phi$ is a critical point of the length function 
$L_m(\phi)$. 

We conclude that for each $1 \le m \le  n/2  $,
the billiard map on a $n$-symmetric oval $\Gamma$ has at least two geometrically distinct trajectories of period
$\tilde n$ corresponding to the maximum and minimum of $L_m$ in $[0,\frac{2\pi}{n}]$. 
As a matter of fact, except for $\tilde n = 2$, each $\tilde n$-polygon gives rise to two different trajectories by reversing the orientation.  
Making the distinction between the orientations,  we can extend the construction of the regular polygons to all $1 \le m \le n-1$ by defining $P_{n-m} (\phi) = P_{m}(\phi)$.
Moreover, if $\tilde n \ne n$, we have $\gdc(m,n) >1$ distinct polygons. In summary, we have

\begin{prop}
\label{prop:traj}
Let $\Gamma \in \mathcal{N}$ be a $n$-symmetric  oval.
Given $1\leq m\leq n-1$,
there are at least $2\gcd (n,m)$ trajectories of period $ \tilde n = n/\gcd(n,m)$. 
Half of them {correspond to regular $\tilde n$-polygons} of maximum perimeter and the other half to minimum perimeter.
The angle between an arbitrary segment of the trajectory and the tangent vector at the vertex is $\alpha_m= \frac{m\pi}{n}$.
\end{prop}

On the other hand, the reflexion law implies that if a regular polygon $P_m(\phi)$ corresponds to a trajectory,  
$\Gamma(\phi)$ is parallel to the normal vector
$\eta (\phi)$.
So from Equation~\ref{eq:curva em g}  the symmetric periodic trajectories correspond to the critical points 
of the support function. 
Then from each point $\Gamma(\phi_0)$ such that $g'(\phi_0) = 0$ there are periodic trajectories 
for each $1 \le m \le n-1$.
By symmetry $g'(\phi_0 + \frac{2k\pi}{n}) = 0$  and so the points $\Gamma(\phi_0 + \frac{2k\pi}{n})$ also belong to periodic trajectories. In fact, all regular polygons with vertices on these points correspond to periodic trajectories.
The polygonal trajectories share the same circumscribed circle of radius $g(\phi_0)$ and are all determined by $\phi_0$.

\begin{defin}
For $1\le m \le n-1$ 
let $\phi_0$ be a critical point of the support function $g$ in $S_n \equiv [0,\frac{2\pi}{n})$.
The family of $n$-symmetric orbits associated to $\phi_0$ is 
$$
\left \{
\mathcal{O}(\phi_k, p_m),
\hbox{ for } {m=1, \ldots, n-1} 
\hbox{ and } {k=0, \ldots, \gdc(n,m)-1} 
\right \}
$$  
where the orbit 
$\mathcal{O} (\phi_k = \phi_0 + \frac{k\pi}{n} , p_m = \cos  \frac{m\pi}{n} ) $
goes through the vertex of the  $\tilde n$-polygon  $P_m(\phi_k)$
$$
\left \{ \Gamma(\phi_k),  \Gamma(\phi_k +  \frac{2m\pi}{n}), \ldots,    \Gamma(\phi_k + \frac{2(\tilde{n} -1)m\pi}{n}), 
\Gamma(\phi_k + \frac{2 \tilde{n}m\pi} {n})  = \Gamma(\phi_k) \right\}
$$
\end{defin}

We can restate Proposition \ref{prop:traj} as
\begin{prop}
{Any  $n$-symmetric  oval in $\mathcal{N}$ has at least two families of $n$-symmetric orbits located at the maximum and minimum of the support function and corresponding respectively  to polygons of maximum and minimum perimeter.}
\end{prop}

We will show in the next section that the elements of a family of symmetric orbits  share some  important dynamical properties.

 % poligonos
\section{Linear stability of the symmetric orbits}
\label{sec:linear}

\begin{figure}
\begin{center}
\includegraphics[width=0.3\hsize]{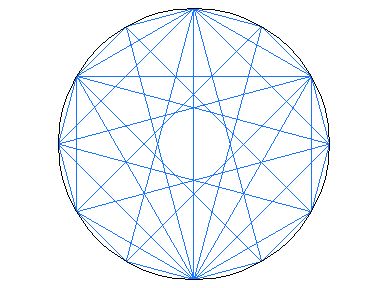}
\hskip 1cm
\includegraphics[width=0.35\hsize]{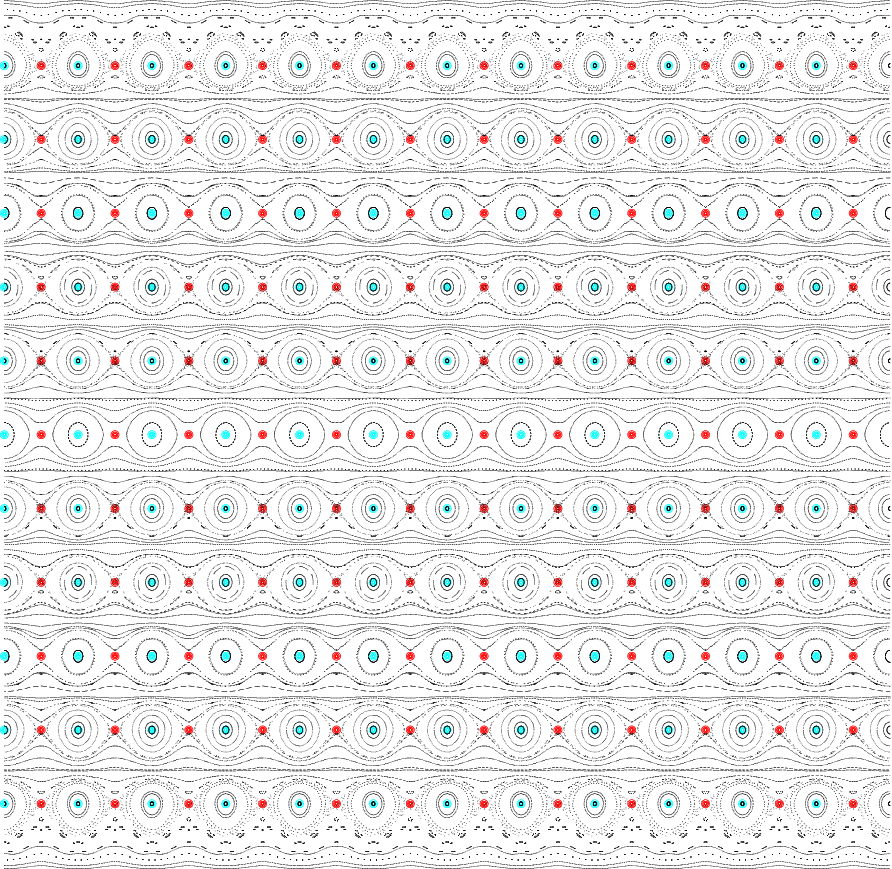} %,valign=c
\end{center}
\caption{Periodic trajectories in a 12-symmetric oval:
%(hyperbolic orbits are red and elliptic blue).
the orbits corresponding to these trajectories  start at the vertical line $\phi =0$. All points $(\frac{2k\pi}{n},  \cos \frac{m\pi}{n})$ for $k=0,\ldots,n-1$ and $m=1,\ldots,n-1$ (blue dots) are periodic and belong to the orbits of the same family. The red dots correspond to the family associated to $\phi=\frac{\pi}{n}$.
The phase space is composed of twelve identical vertical strips, each one contains exactly one blue and one red dot.}
\label{fig:poligonos12}
\end{figure}

In order to investigate the dynamical properties of the symmetric orbits, we will explore how the symmetry of the curve reflects itself on the phase space.
It is clear that the invariance of the curve under the rotation of $\frac{2\pi}{n}$ implies that the billiard map is
invariant under the translation $\phi \to \phi + \frac{2\pi}{n}$ and
thus the phase space consists of $n$ equivalent vertical strips which are clearly seen on 
Figure~\ref{fig:poligonos12}.
Observe also the horizontal symmetry around $p=0$, due to the reversibility.

Given $m$ such that $1 \le m \le n-1$,  we define the equivalence relation
$$ (\phi_1,p_1) \sim_m (\phi_2,p_2) \Leftrightarrow \phi_2 = \phi_1 \hbox{ mod } \frac{2m\pi}{n}$$  
This induces a map $T_m$ on the cylinder 
%${\obs \Omega_m} =  \mathbb{R}/{ {\frac{2m\pi}{n}} \mathbb{Z}} \times (-1,1)$
${ \Omega_m} = S_n \times (-1,1)$
given by
$
T_{m}([(\phi,p)]_m) = [T(\phi,p)]_m$.
This map is also a  diffeomorfism with the quotient topology and we will call it {\em the $m-$quotient map}. 
We will sometimes use simply $(\phi,p)$ to also  denote the 
equivalence class $[(\phi,p)]_m$, when applicable.

By construction, a point $(\phi_{0}, p _m)$ with  $p_m = \cos  {\frac{m\pi}{n}}$ belongs to a periodic orbit associated to
a polygon  $P_{m}$ if and only if $[( \phi_{0}, p_m)]_m $ is a fixed point of $T_{m}$.
The result of the previous section also holds in this context.
More specifically, 
given a $n$-symmetric oval $\Gamma $ for every $1\leq m\leq n-1$, 
the fixed points of the quotient map $T_m$  correspond to the critical points of the support function $g$ and so we have at least two fixed points.

%%%%%%%%%%%%%%%%%%%%%%%%%%%%%%%%%%%%%%%%%%%%%
\begin{comment}
Its derivative is given by \cite{xxx}
\begin{equation}
\label{eq: derivada do bilhar (phi,p)}
DT_{(\phi_0,p_0)} = 
\left(
\begin{array}{cc}
\displaystyle
\frac{L -R_0 \sin \alpha_0}{R_0 \sin \alpha_1}
&
\displaystyle
\frac{-L}{R_0 \sin \alpha_0 \sin \alpha_1 }    \\
& 
\\
\displaystyle
\frac{-(L - R_0 \sin \alpha_0 - R_1 \sin \alpha_1)} {R_0} 
&
\displaystyle
\frac{L - R_1 \sin \alpha_1 }{R_1 \sin \alpha_0 }
\end{array}
\right)
\end{equation}
where $R_0 = R(\phi_0)$ and $L=L(\phi_0, \phi_1)$.
The measure $R(\phi)  d\phi \, dp$ is $T$-invariant. Clearly, the billiard map can also be given using the variable $s$ instead of $\phi$, we will denote this map by $T(s,p)$. In this case, the Lebesgue measure $ds\,dp$ is preserved. 
\end{comment}
%%%%%%%%%%%%%%%%%%%%%%%%%%%%%%%%%%%%%%%%%%%%%%%%%

\begin{prop}
Consider an oval $\Gamma \in \mathcal{N}$ and  let $\phi_0$ be a critical point of its support function $g \in {\mathcal{G}}$. 
Then the orbits of the family associated to $\phi_0$ are elliptic (hyperbolic) if and only if $\phi_0$ is a non degenerate minimum (maximum) of $g$.  
\label{prop:estabilidade linear}
\end{prop}

\begin{prova}
If $\phi_0$ is a critical point of the support function, $(\phi_0, p_m = \cos \frac{m\pi}{n})$ is a fixed point of $T_m$. 
It is clear that
$DT_m ( \phi,p) = DT (\phi,p) $. 
The expression of $DT$ is well known \cite{markarian} and at the fixed point $(\phi_0, p_m)$ it is given by
\begin{equation}
\label{eq: derivada do bilhar (phi,p)}
DT_m{(\phi_0,p_m)} = DT{(\phi_0,p_m)} =
\left(
\begin{array}{cc}
\displaystyle
\frac{L_m(\phi_0) -R_0 \sin \alpha_m} {R_0  \sin  \alpha_m}
&
\displaystyle
\frac{-L_m (\phi_0)}{R (\phi_0) \sin^2 \alpha_m  }    \\
& 
\\
\displaystyle
\frac{-(L_m(\phi_0) - 2 R_0 \sin \alpha_m } {R_0} 
&
\displaystyle
\frac{L_m (\phi_0) - R_0 \sin \alpha_m}{R_0 \sin \alpha_m }
\end{array}
\right)
\end{equation}
where $R_0= R(\phi_0)$, $\alpha_m = \frac{m\pi}{n}$ and  $ L_m (\phi) =|| \Gamma(\phi+ 2 \alpha_m ) - \Gamma(\phi) ||$. Clearly we have
\begin{equation}
\det  DT_m( \phi_{0}, p_m)  =1
\hbox{ and }
\hbox{ Tr } DT_m( \phi_{0}, p_m)  
=\frac{2L_m( \phi_{0}) }{R( \phi_{0}) \sin \alpha_m}-2
\label{eq:tracoDT}
\end{equation}
and so the point is
\begin{itemize}
 \item  elliptic if and only if 
$L_m\left( \phi_{0}\right)
< 2 R_0 \sin {\frac{m\pi}{n}} $
 \item parabolic if and only if
$L_m(\phi_{0}) =2R\left( \phi_{0}\right) \sin  {\frac{m\pi}{n}} $
\item  hyperbolic if and only if
$L_m\left( \phi_{0}\right) >
2 R_0 \sin {\frac{m\pi}{n}} $ 
\end{itemize}

So  the point $(\phi_0, p_m)$ is elliptic if and only if $g(\phi_0) < R_0$, hyperbolic if and only $g(\phi_0) > R_0$ and parabolic if and only if $g(\phi_0) = R_0$.   In particular, this implies that the linear stability of all the $n$-symmetric associated to $\phi_0$ depends only on $\phi_0$ and not on $m$. 
As $R_0 = g(\phi_0) + g''(\phi_0)$ we have that the set of $n$-symmetric orbits associated to $\phi_0$ is hyperbolic (elliptic) if and only if $\phi_0$ is a non degenerate maximum (minimum) of the support function. 
\end{prova}

Now as  $\Gamma(\phi) = g(\phi) \eta(\phi) - g'(\phi) \sigma(\phi)$, we have that
$$
L_m^2(\phi) %=2 (g(\phi)^2+g'(\phi)^2) \left (1- \cos \frac{2m\pi}{n}\right)
= 4   (g(\phi)^2+g'(\phi)^2)  \sin^2 {\frac{m\pi}{n}} 
$$
and at a critical point  $\phi_{0}$ 
$$L_m(\phi_0) = 2 g(\phi_0) \sin  {\frac{m\pi}{n}}$$

Strictly speaking, if the eigenvalue of the fixed point is an $m$-root of the unity, the orbit itself would be parabolic. 
Nevertheless, in this case, we will say that the orbit is elliptical because, as we will see in the next sections, its stability is completely determined by the stability of the fixed point.

Clearly, the subset of Morse support functions  is an open and dense subset of the support functions $\mathcal{G}$ with at least two sets of symmetric periodic orbits, one composed of linearly  elliptic orbits and the other one of hyperbolic orbits, there are no parabolic symmetric orbits in this case.
Nevertheless it is of course possible that a symmetrical billiard has an elliptical or hyperbolic family even if it is not Morse, that is, its support function may also have degenerate critical points, isolated or not.
Anyway, having a non degenerate minimum (and/or a degenerate maximum) is clearly an open property in the set of support functions $\mathcal{G}$. On the other hand, the density of Morse functions implies that this is a dense property.
Thus we have

\begin{prop}
\label{prop:2familias}
The set of $n$-symmetric ovals  
that have at least one family of elliptic symmetric orbits and one family of hyperbolic symmetric orbits 
is
open and dense in $\mathcal{N}$
with the topology induced by $\mathcal{G}$.
\label{prop:2familias}
\end{prop}

It is worthwhile to note that the results in this section hold for $C^2$ support functions. We also stress that, for symmetric billiards,  our result implies the existence of elliptic orbits, which are not guaranteed by Birkoff's Min-Max.

%%%%%%%%%%%%%%%%%%%%%%%%%%%%%%%%%%%%%%%%%%%%%%%%%%%%%%%%%%
 %linear
\section{Stability of elliptic orbits} \label{sec:elip}

In this section we will analyze the stability of the symmetric elliptical orbits. 
This  stability can be established  by Birkhoff's Normal Form \cite{birk} and Moser's Twist Theorem \cite{moser}, which states conditions for the existence of infinitely many invariant curves surrounding the points of the orbit resulting in elliptic islands of positive Lebesgue measure (Moser stability). In this approach, we need to switch to canonical variables $(s,p)$, where $s$ is the arclenght parameter of the boundary $\Gamma$, i.e. $ds=R \,d\phi$. We will interplay between the variables $s$ and 
$\phi$ and in either case denote the curve as $\Gamma$ and the billiard map as $T$.

Given  $1 \le m \le  n-1$, let $(\phi_0 = \phi(s_0), p_m=\cos \dfrac{m}{n}\pi )$ be an elliptic fixed point of the quotient map $T_m$ associated to an oval $\Gamma$.
Using canonical variables $(s,p)$ and the coordinates of the Jordan Form of  
$DT_m(s_0,p_m)$ at the fixed point, the complex form of the map is written as \cite{moeckel}
\begin{equation}
z\rightarrow \lambda \left( z+c_{20}z^{2}+c_{11}z\bar{z}+c_{02}\bar{z}%
^{2}+c_{30}z^{3}+c_{21}z^{2}\bar{z}+c_{12}z\bar{z}^{2}+c_{03}\bar{z}%
^{3}\right) +O\left( \left\vert z\right\vert ^{4}\right)
\label{eq: primeira complexificao}
\end{equation}
where $\lambda =\cos\zeta \pm i\sin \zeta $ are the eigenvalues of $DT_m$. 
If the fixed point is non resonant, by the 
Birkhoff's Normal Form Theorem, there is a $C^{\infty }$ coordinate change bringing the map into the form
$$
z\rightarrow e^{i\left( \zeta +\tau\left\vert z\right\vert ^{2}\right)
}z+O\left( \left\vert z\right\vert ^{4}\right) .
$$
Moser's Twist Theorem then assures that if the 
first Birkoff Coefficient (Twist Coefficient) $\tau$ is non zero, the fixed point is stable. 

%definicao ressonante%
\begin{defin}
An elliptic fixed point  of a conservative bi-dimensional map 
is resonant (of order 4) if its eigenvalues 
satisfy $\lambda^j = 1$ for $j = 1, 2, 3$ or $4$.
Otherwise we say that the fixed point is non resonant.
\end{defin}

From the proof of Proposition \ref{prop:estabilidade linear} we have that 
$\displaystyle
\cos \zeta = \frac{L_m(\phi_0)}{R_0 \sin \alpha_m} -1
$.
%where  $R_0 = R(\phi_0)$, $\alpha_m = \frac{m}{n} \pi$and $ L_m (\phi_0) =|| \Gamma(\phi_0+ 2 \pi \frac{m}{n}) - \Gamma(\phi_0) ||$.
If the point is elliptic,  $0 < L_m(\phi_0)< 2 R_0 \sin \alpha_m$ and $\lambda, \lambda^2 \ne  1$.
So the only resonance conditions are 
\begin{eqnarray*}
\lambda^3 = 1 
& \Leftrightarrow 
\ \zeta = \pm 2\pi/3   
& \Leftrightarrow
 \ 2 L_m(\phi_0) = R_0 \sin \alpha_m
\\
\lambda^4 = 1 
&\Leftrightarrow 
\ \zeta = \pm \pi/2   
& \Leftrightarrow \  L_m(\phi_0) =  R_0 \sin \alpha_m
\end{eqnarray*}
As $L_m(\phi_0)=2 g(\phi_{0}) \sin \alpha_m $ and
$R_0 =g(\phi_{0}) +g^{\prime \prime }(\phi_{0})$ 
we can rewrite the resonance conditions in terms of the support function
\begin{equation}
\lambda^3 = 1 \Leftrightarrow 3 g (\phi_{0}) - g^{\prime \prime }(\phi_{0}) = 0 
\ \ , \ \
\lambda^4 = 1 \Leftrightarrow g(\phi_{0}) -g^{\prime \prime }(\phi_{0}) = 0 \label{eqn:g-resonante}
\end{equation}
We say that $\phi_0$ is a resonant minimum of $g$ if at least one of conditions above holds,  otherwise, we say that we have a non resonant minimum. 
As these equations do not  depend on $m$, if $\phi_0$ is a non-resonant minimum, the fixed points $(\phi_0, \alpha_m)$ will be non-resonant for all $m$.
However, as the eigenvalues of the orbit of $(\phi_0,\alpha_m)$ are $\lambda^m$, the orbit itself may be resonant even if the fixed point is not.
Nevertheless we define
\begin{defin}
A family of elliptic symmetric orbits  is resonant if it is associated to a resonant minimum of the support function. Otherwise we say that the family is non resonant. 
\end{defin}

%%%%%%%%%%%%%%%%%%%%%%%%%%%%%%%%%%%%%%%%%%%%%%%%%%%%%
%lema: g morse com ressonancia pode ser perturbada

\begin{lem} \label{lem:res}
 If $g \in {\mathcal{G}}$ 
 is a support function with one or more 
resonant minima, there is $g_{\varepsilon}$ also in ${\mathcal{G}}$ and {$C^2$-close} 
to $g$ such that all its non degenerate minima are non resonant and coincide with the non degenerate minima of $g$ .
\end{lem}

\begin{prova} 
Let $g_{\varepsilon} = {g} + \varepsilon$.
Obvioulsy, $g^{\prime}_{\varepsilon }= {g}^{\prime}$ and
 $g^{\prime\prime}_{\varepsilon }= {g}^{\prime\prime}$  for all $\phi$ and so they have the same minimum (and maximum) points. Moreover
$\| {g} - g_{\varepsilon} \| = |\varepsilon| $.
It is clear that we can choose $\varepsilon \sim 0$ such that the non degenerate minima of $g_{\varepsilon}$ are all non resonant, i.e.
$g_{\varepsilon} \ne g_{\varepsilon}^{\prime\prime}$,  $g_{\varepsilon} \ne  g_{\varepsilon}^{\prime\prime}/3$ 
simultaneously at all non degenerate minimum points. 
\end{prova}

At this point it is worthwhile to observe that adding a constant to the support function results in a normal perturbation of the associated oval. 

%%%%%%%%%%%%%%%%%%%%%%%%%%%%%%%%%%%%%%%%%%%%%

\begin{prop}
The set of ovals such that all families of symmetric elliptic orbits are non resonant is open and dense in 
$\mathcal{N}$
\label{prop: nao ressonantes aberto e denso}
\end{prop}

\begin{prova}%[Proof of Proposition \ref{prop: nao ressonantes aberto e denso}]
Let $g$ be a support function with a non resonant minimum at a point $\phi_0$, i.e.
$g'(\phi_0)=0$, $g''(\phi_0) >0$ , $g''(\phi_0)\ne g(\phi_0)$ and $g''(\phi_0)\ne 3 g(\phi_0)$.  
Then any support function 
$\tilde{g}\in \mathcal{G} $ 
with $\left\vert \left\vert \tilde{g}-g\right\vert \right\vert_{2}<\delta $ 
will have a non resonant minimum $\widetilde{\phi}_0 \sim \phi_0$ for $\delta$ small enough
{and so this is an open property}.
Since non degenerate points are isolated and finite, the argument holds for all non degenerate points.

Given an oval $\Gamma$ with support function $g$  and $\delta > 0$, 
we can choose a Morse function $\tilde{g} \in \tilde{\mathcal{G}}$ 
such that
$\left\vert \left\vert {g}- \tilde{g}\right\vert \right\vert_{2}<\delta /2$.
Now, if necessary,  let $g_{\varepsilon} = \tilde{g} + \varepsilon$, with $|\varepsilon|<\delta/2$ chosen as in 
Lemma~\ref{lem:res} above  and $\Gamma_{\varepsilon}$ the associated curve.
We have
$\| \Gamma_{\varepsilon }-\Gamma \|_{2}=
\| g_{\varepsilon }-g \|_{2}<
\| g- \tilde{g}\|_{2}
+
\| \tilde{g}-g_{\varepsilon } \|_{2} <
\delta $ which proves that {the subset of ovals such that all families of elliptical symmetric orbits are non resonant} is dense in $\mathcal{N}$. 
\end{prova}

%%%%%%%%%%%%%%%%%%%%%%show is open%%%%%%%%%%%%%%%%
%Calculo do tau 1
%%%%%%%%%%%%%%%%%%%%%%%%%%%%%%%%%%%%%%%%%%%%%%%

Once we have a non resonant elliptic fixed point we can go on and determine
the Twist Coefficient which depends on the derivatives of the map up to order 3 \cite{moeckel,ss}
\begin{equation}
\tau=\mbox{Im} ( c_{21}) 
+\frac{\sin \zeta }{\cos \zeta-1}
\left( 
3 \vert c_{20} \vert ^{2}
+\frac{2\cos \zeta -1}{2\cos\zeta +1} \vert c_{02}\vert ^{2}
\right)
\label{eq: formula coef de birkh}
\end{equation}

Obtaining the 3-jet of 
$T_m $ at $(s_{0} = s(\phi_0), p_m= \cos \alpha_m)$ and the $c_{ij}$'s is a straightforward but long calculation which, using  the symbolic computation software  Maple,  leads to
\begin{eqnarray}
\tau &=&
-\frac{1}{8 R_0 \sin ^{3}\alpha_m }
+\frac{3 \cos^{2} \alpha_m }{8\sin^{2} \alpha_0 \left( 2L_m(\phi_0) -R_0 \sin \alpha_m \right) }
 \label{eq:tau1} \\
%&&-\frac{L_m(\phi_0) \left( 7 L_m(\phi_0) -4 R_0 \sin \alpha_m \right) }{8 \left( L_m(\phi_0) -2 R_0 \sin \alpha_m \right)^{2} \left( 2 L_m(\phi_0) -R_0 \sin \alpha_m \right)} {R^{\prime}_0}^2 \nonumber \\
&&-\frac{1}{8}\frac{L_m(\phi_0)}{(L_m(\phi_0)-2 R_0 \sin \alpha_m)^2} 
\left(
3 + \frac{L_m(\phi_0)-R_0 \sin \alpha_m }{2 L_m(\phi_0)-R_0 \sin \alpha_m }
\right)
{R^{\prime}_0}^2 \nonumber \\
&&-\frac{L_m(\phi_0) }{8\sin \alpha_m \left(
L_m(\phi_0) -2R_0 \sin \alpha_m\right) }R^{\prime \prime}_0
\nonumber
\end{eqnarray}%
where $R'_0 = \frac{dR}{ds}(s_0)$ , $R''_0 = \frac{d^2R}{ds^2}(s_0)$ and,  as before, $\alpha_m = \frac{m\pi}{n}$ and $R_0 = R(\phi_0)$.

%%%%%%%%%%%%%%%%%%%%%%%%%%%%%%%%%%%%%%
% no maximo um periodo com tau 1 nulo
%%%%%%%%%%%%%%%%%%%%%%%%%%%%%%%%%%

\begin{prop} 
Let $\Gamma \in {\mathcal N}$
 be an oval with a family of elliptical non resonant symmetric orbits. 
The twist coefficient of these orbits can only be zero for at most one period, i.e. for  the orbits associated to at most one value of $1\le m \le  n/2 $.
\label{prop:max1tau0}
\end{prop}

\begin{prova}
Let the family be associated to $\phi_0$. 
As $L_m(\phi_0) = 2 g (\phi_0) \sin \alpha_m$, the first Birkoff coefficient associated to an orbit in the family, given by \ref{eq:tau1} depends on $m$ only through  
$\alpha_m = \frac{m\pi}{n}$.
\begin{eqnarray*}
\tau&=&
-\frac{1}{8 R_0 \sin ^{3}\alpha_m }
+\frac{3 \cos^{2} \alpha_m }{8\sin^{3} \alpha_m ( 4 g_0 -R_0) }
 \\
&&-\frac{g_0}{16 \sin \alpha_m (g_0-R_0 )^2} 
\left(
3 + \frac{2g_0-R_0  }{4 g_0-R_0 }
\right)
{R^{\prime}_0}^2  \\
&&-\frac{ g_0 }{8\sin \alpha_m (g_0 -R_0) }R^{\prime \prime}_0
\end{eqnarray*}
where $g_0 = g(\phi_0)$.

Then  $\tau=0$ if and only if
\begin{equation}
4\frac{g_0-R_0 }{R_0 ( 4 g_0 -R_0) }
+\frac{3 \sin^{2} \alpha_m }{( 4 g_0 -R_0) }
 +\frac{g_0 \sin ^{2}\alpha_m }{2  (g_0-R_0 )^2} 
\left(
3 + \frac{2g_0-R_0  }{4 g_0-R_0 }\right){R^{\prime}_0}^2 
+\frac{ g_0 \sin ^{2}\alpha_m }{ (g_0 -R_0) }R^{\prime \prime}_0
=0
\label{eq:tau0linear}
\end{equation}
which is a linear equation in $\sin^2 \alpha_m$. 
Clearly, if $\phi_0$ (and so $g_0$ and $R_0$) is given and as $1\le m\le  n/2 $, this equation has at most one solution.
\end{prova}

An immediate consequence of this proposition is that having a symmetric Moser stable orbit  is an open an dense property for ovals with symmetry $n \ge 4$. 

\begin{comment}
-- ter ilha eh denso porque as Morse  nao resonantes tem ilha 
-- ter ilha eh aberto eh geral demais. o que sabemos mostrar que eh aberto eh ter uma orbita simetrica nao ressonante com twist nao nulo 
\end{comment}

\begin{thm}
\label{thm:teorema1}
There is an open and dense subset of $C^{\infty}$ ovals with symmetry $n\ge 4$ with an elliptical island. 
\end{thm}

\begin{prova} 
By Proposition \ref{prop:max1tau0} above, if the symmetry $n$ is at least $4$ it is enough to find an open and dense  set  of ovals with a non resonant elliptic family since at least one of the associated symmetric orbits will have a non zero twist coefficient and so will be Moser stable. Proposition \ref{prop: nao ressonantes aberto e denso} gives us such a set.
\end{prova}

With regard to the restriction of $n \ge 4$ in the theorem above, we observe
 if $n=2$ the only symmetric orbits have period two and  
if $n=3$  we only have symmetric orbits of period three. 
% n impar = largura constante (as orbitas de periodo 2 nao sao simetricas...)
In both case, we can have a zero twist coefficient. The generic stability of period 2 orbits has been established in \cite{ilhas,ss} for general curves, and the result also holds within symmetric curves. Anyway Theorem~\ref{thm:teorema2} bellow applies in these cases.

%%%%%%%%%%%%%%%%%%%%%%%%%%%%%%%%%%%%%%%%%%%%%
%Perturbando o coeficiente de twist
%%%%%%%%%%%%%%%%%%%%%%%%%%%%%%%%%%%%%%%

{
Finally, in order to prove that a zero twist coefficient is a rather rare fact, we can use the same argument of 
\cite{ilhas}, which relies on the remark that $\tau$ depends  linearly on $R''_0$
and its  coefficient
$$-\frac{L_0}{8\sin \alpha_m \left(
L_0 -2R_0 \sin \alpha_m
\right) }\neq 0$$
This implies that a small perturbation of $R^{\prime\prime}$  is enough to slightly modify this term in order to obtain a non-zero twist coefficient

}

\begin{lem}
\label{prop:aprox oval tau dif 0}
Let  $\Gamma \in {\mathcal{N}}$ be
such that the  $(\phi_0,\alpha_m)$ corresponds to a  non resonant elliptic fixed point of the quotient map $T_m$  with a zero first Birkhoff coefficient $\tau$.
Then there is an oval  $\Gamma_{\epsilon} \in \mathcal{N}$
arbitrarily close to $\Gamma$,  such that the point $(\phi_0, \alpha_m)$ also corresponds to a non resonant elliptic point of the map ${T}_m$,  but its twist coefficient is non zero.
\end{lem}
\begin{prova} Let $g$ be the support function of $\Gamma$.
We can assume, without any loss of generality, that the elliptic family with a zero twist coefficient is associated to $\phi_0 = 0$,
which is an isolated minimum,  
and so $g'(0) = 0$ and $g^{\prime\prime} > 0 $. 

Given $0<\delta_1 < \delta_2$  let $\rho(x)$ be a $C^{\infty}$ 
function with compact support such that:
\begin{lista}
\item
$\rho(-x) = \rho(x)$
\item
$\rho (x) = 0$ for $|x| > \delta_2$
\item
$\rho (x) = 1$ for $|x| < \delta_1$ (all the derivatives of $\rho$ vanish in a neighborhood of $x=0$)
\item
$ 0<-\rho^{\prime}(x) < C$ for $\delta_1 < x < \delta_2$
%\item nao sei se eh morse, mas acho que tanto faz...
\end{lista}
There are well know functions with these properties which are frequently used in the construction of 
partitions of the unity (smooth indicator or bump functions).

We choose $\delta_2$ such that the support of $\rho$ does not contain any  other critical points of $g$, 
i.e., $g'(\phi) \ne 0$ for  $0<|\phi|<\delta_2$.
Now let   $g_{\epsilon}(\phi) = g(\phi) +  \epsilon \phi^4 \rho(\phi)$ with $\epsilon \ne 0$ small. Then
\begin{lista}
\item
 $g_{\epsilon}(0) = g (0)$ 
\item 
 $g^{\prime}_{\epsilon}(0) = g^{\prime} (0) = 0$, 
$g^{\prime \prime}_{\epsilon}(0) = g^{\prime \prime} (0)$, 
$g^{\prime \prime \prime}_{\epsilon}(0) = g ^{\prime \prime \prime}(0)$
\item
 $g^{\prime}_{\epsilon}(\phi) = g^{\prime} (\phi )$ for $ 0 \le |\phi| \le \delta_1$ (and for $|\phi| >\delta_2$).
\item
 $g^{\prime}_{\epsilon}(\phi) = g^{\prime} (\phi ) + \epsilon \phi^3 (4 \rho(\phi) + \phi \rho^{\prime}(\phi) )$ for
$ \delta_1 < |\phi|  < \delta_2$. As $g^{\prime} \ne 0 $, $0 < \rho < 1$ and $\rho^{\prime}$ is bounded in this region, we can choose $\epsilon$ small enough in order that 
$g^{\prime}_\epsilon$ is also not zero and  $g_{\epsilon}$ and $g$ have the same critical points in $M =  [0,2\pi/n)$.
\item
$g^{\prime \prime \prime \prime}_{\epsilon}(0) = g ^{\prime \prime \prime \prime}(0) + 24 \epsilon $
\end{lista}

Let $\Gamma_{\epsilon}$ be the oval in $\mathcal{N}$ which corresponds to the 
support function $g_{\epsilon}$.
As ${\Gamma}_{\epsilon}$ and $\Gamma$ also have a third order contact at $\phi_0$, 
the point $ ( \phi_{0},\alpha_m)$
is also a non resonant elliptic fixed point of the map quotient map associated to
$\Gamma_{\epsilon}$. 
Moreover   
$$ || {\Gamma}_{\epsilon}(\phi_0 + \frac{2\pi m}{n}) - {\Gamma}_{\epsilon}(\phi_0) ||
=   || {\Gamma}(\phi_0 + \frac{2\pi m}{n}) - {\Gamma}(\phi_0) ||
= L_m(\phi_0) = L_0$$
and ${R}_{\epsilon}(\phi_0) = {R}(\phi_0)  = R_0 $ as well as the first derivatives.
 
If ${\tau}_{\epsilon}$ and  ${\tau}$ denote the first coefficient of the 
non resonant elliptic fixed point $(\phi_0, p_m)$ for the quotient map  $T_m$ associated respectively to $\Gamma_{\epsilon}$ and $\Gamma$
we have
$$
{\tau}_{\epsilon} -\tau = 
- \frac{L_0 }{
8 R_0^2 \sin \alpha_m ( L_0 -2R_0 \sin \alpha_m) }%
\left( 
\frac {d^2 {R}_{\epsilon}}{d{\phi}^2}(\phi_0) - \frac {d^2 {R}}{d{\phi}^2}(\phi_0)
\right)
$$
As by hypothesis  $\tau=0$ and we have $\frac{d^2 R}{d\phi^2} = g'' + g''''$, by choosing ${g}_{\epsilon}$ such that 
${g}'''''_{\epsilon} (\phi_0)  \ne {g}''''' (\phi_0) $
we obtain 
${\tau}_{\epsilon}\neq 0$ which concludes the proof.
Moreover, as changing $g''''$ wont change the independent and the  $(R')^2$ terms in $\tau$, we can do this without causing the twist coefficients for the other values of $m$ to vanish.
\end{prova}

We conclude then that the set of Morse $n$-symmetric ovals, such that all symmetric orbits are non resonant and have a non zero twist coefficient is dense in $\tilde{\mathcal N}$ and so also in $\mathcal N$. As these properties are clearly open we have
%%%%%%%%%%%%%%%%%%%%%%%%%%%%%%%%%%%%%%%%%%%%%%%%%%%
%%%%%%%%%%%%%%%%%%%%%%%%%%%%%%%%%%%%%%%%%%%%%%%%%%%
%{ter todas as orbitas simetricas elipticas Moser estaveis eh uma propriedade aberta e densa entre as ovais simetricas? } 
\begin{thm}
\label{thm:teorema2}
There is an open and dense subset of the $n$-symmetric ovals  where all the elliptic symmetric orbits are Moser stable.
\end{thm}

\begin{comment}
 In fact, our approach using Birkoff's Normal Form and Moser's Twist Theorem works for $C^{4}$ area preserving maps and so theorems \ref{thm:teorema1} and \ref{thm:teorema2} hold for $C^5$ curves.
\notaIL{Syok}{An alternative approach uses Herman's Last Geometric Theorem \cite{fayadkrikorian} and holds for $C^{\infty}$ maps with a diophantine fixed point.}
\end{comment}

 %elip
%5-hyp
\section{Hyperbolic Orbits} \label{sec:hyp}

Proposition \ref{prop:2familias} guarantees the existence of an open and dense set of ovals with at least one family of hyperbolic orbits. We will show that Theorem~2 of Dias Carneiro et al. \cite{periodicas}
which states that  generically, for billiards on ovals, the stable and unstable manifolds of two
hyperbolic points either do not meet or have at least one transversal homo/heteroclinic
connection holds in the context of symmetric curves.

\begin{prop} Consider a symmetric curve $\Gamma$ and
let $(\phi_0, \alpha_0)$ and $(\phi_1,\alpha_1)$ be two hyperbolic points and let us assume that there is a tangent 
hetero(homo)clinic point
$(\phi_*,\alpha_*) \in W^s (\phi_0, \alpha_0) \cap W^u(\phi_1,\alpha_1)$.    
Then there is a symmetric curve $\Gamma_{\epsilon}$ close to $\Gamma$ such that  $(\phi_*,\alpha_*)$ is a transversal hetero(homo)clinic intersection.
\end{prop}

\begin{prova}
The original  proof uses a normal perturbation of the curve.
We can repeat the construction in our case using a local perturbation of the support function.

As $(\phi_*,\alpha_*) \in W^s (\phi_0, \alpha_0)$, there is a neighborhood of this point that does not contain any other points in its forward orbit, and the same is true for the backward orbit. 

We can then find an interval $I \ni \phi_*$  such that if 
$(\phi,\alpha)$ is in the orbit of  $(\phi_*,\alpha_*)$, then $\phi \notin I$. Moreover, the twist property ensures that we can assume that the tangency is not vertical and so $I$ can be chosen small enough that $W^s$ and $W^u$ are local graphs given by $\alpha = \alpha^u(\phi)$ and $\alpha=\alpha^s(\phi)$ with $\alpha^u(\phi_*) =   \alpha^s(\phi_*) = \alpha_*$ and
$\frac{d\alpha^u}{d\phi} (\phi_*) = \frac{d\alpha^s}{d\phi} (\phi_*) = \alpha_*'$ .

To the single tangent vector $(1,\alpha_*')$ corresponds a pencil of rays focusing forward and backward respectively at distances
$$
d_+= \frac{R_* \sin \alpha_*}{1+\alpha_*'} \ \  \hbox{ and }  \ \ d_-= \frac{R_* \sin \alpha_*}{1-\alpha_*'}
$$
where $R_* = R(\phi_*)$ is the curvature radius at $\phi_*$.

We observe that $\phi_*$ does is not necessarily a symmetric periodic point, nor are $\phi_0$ and $\phi_1$, however by  symmetry the construction will apply to all their rotations by $\frac{2\pi}{n}$. we observe that we can choose the perturbation in order to stay away of the symmetric orbits.

Without any loss of generality, we assume in what follows that $\phi_* = 0$. 
If $g$ is the support function of $\Gamma$ we consider the curve 
$\Gamma_{\epsilon}$ given by 
$$
g_{\epsilon}(\phi) = g(\phi) + \epsilon \phi^2 \rho(\phi)
$$
where the function $\rho$ is a smooth bump
function as the one used in the proof of Proposition~\ref{prop:aprox oval tau dif 0}.  
We can adjust its support in order to guarantee that no symmetry equivalent of the hyperbolic points $\phi_0$ and $\phi_1$ are inside it and also to match the interval $I$ above.
As $g_{\epsilon}(0) = g(0)$ and $g'_{\epsilon}(0) = g'(0)$,
$\Gamma_\epsilon$ and $\Gamma$ have a first order contact at $\phi_*$. Moreover they coincide outside the small intervals defined by $I$,  so $(\phi_*,\alpha_*)$ is also a hetero(homo)clinic point for the billiard map in $\Gamma_{\epsilon}$ with the same properties. In particular the stable and unstable curves are also local graphs in its neighborhood with 
$\alpha_{\epsilon}^u(\phi_*)=\alpha_{\epsilon}^s (\phi_*)=\alpha_*$
and $\frac{d\alpha_{\epsilon}^u}{d\phi} (\phi_*) =  (\alpha^u_*)'$ , $\frac{d\alpha_{\epsilon}^s}{d\phi} (\phi_*) = (\alpha^s_*)'$ .

On the other hand, as 
$g^{\prime\prime}_{\epsilon} 
= g^{\prime\prime}+ \epsilon (\phi^2 \rho^{\prime\prime} + 4 \phi \rho' + 2\rho)$ 
the curvature radius of $\Gamma_{\epsilon}$ at the intersection $\phi_*$ is 
$$
R_{\epsilon}= g_{\epsilon}(0) + g^{\prime\prime}_{\epsilon} (0)
= R_*+ 2 \epsilon  
$$

The tangent vector $(1,  (\alpha^u_{\epsilon})')$ focus backwards at a distance 
$
(d^u_{\epsilon})_-= \frac{R_{\epsilon} \sin \alpha_*}{1- (\alpha^u_{\epsilon})'}
$
which, as the backward orbit is unchanged, must be 
equal to $d_-$ yielding to 
$$
 (\alpha^u_{\epsilon})' =  \alpha_{*}' - 2\epsilon \, \frac{1-  \alpha_{*}'}{R_*} 
$$
The same argument applies to the forward orbit implying that 
$  (d^s_{\epsilon})_{+} = \frac{R_{\epsilon} \sin \alpha_*}{1- (\alpha^s_{\epsilon})'} = d_+
$ 
which gives
$$
(\alpha^s_{\epsilon})' =  \alpha_{*}' + 2\epsilon \, \frac{1 +  \alpha_{*}'}{R_*} \ne  (\alpha^u_{\epsilon})' 
$$
As tangent vectors have different slope the intersection is transverse.
 \end{prova}

As in \cite{periodicas} we observe that we do not prove that every homo/heteroclinic orbit is transversal. We do
know that generically two invariant stable and unstable manifolds either do not intersect or
have at least one transversal homo/heteroclinic orbit, but there can also be tangent orbits.

On the other hand, Xia and Zhang \cite{xia} have proven that  $C^{\infty}$ generically every hyperbolic periodic point of a billiard admits a homoclinic orbit. 
It is then, a challenging question if their argument can be rephrased  in the context of symmetric curves implying the existence of homoclinic points associated to symmetric orbits. If his is true,
having a hyperbolic orbit with a transverse homoclinic point holds in an open and dense subset of the
$n$-symmetric ovals.  %hyp
\section{Conclusion} \label{sec:final}

\notaIL{Syok}{Para revisao, inclusive as figuras...}

The complete description of a typical phase space of a strictly convex billiard (Birkhoff) is still a challenge. 
With this work, we intent to address this question in the scope of symmetric curves. 

We have shown that symmetric billiards in general mix stable (elliptic islands) and unstable (hyberbolic orbits) behavior.  
It is also  well known that convex billiards have invariant spanning curves \cite{laz} close to the bottom and the top of the phase space. 
The figure bellow illustrates what we expect to be the generic phase space of $n$-symmetric ovals by clearly displaying these elements.

\begin{figure}[h]
	\begin{center}
		\includegraphics[height=4cm]{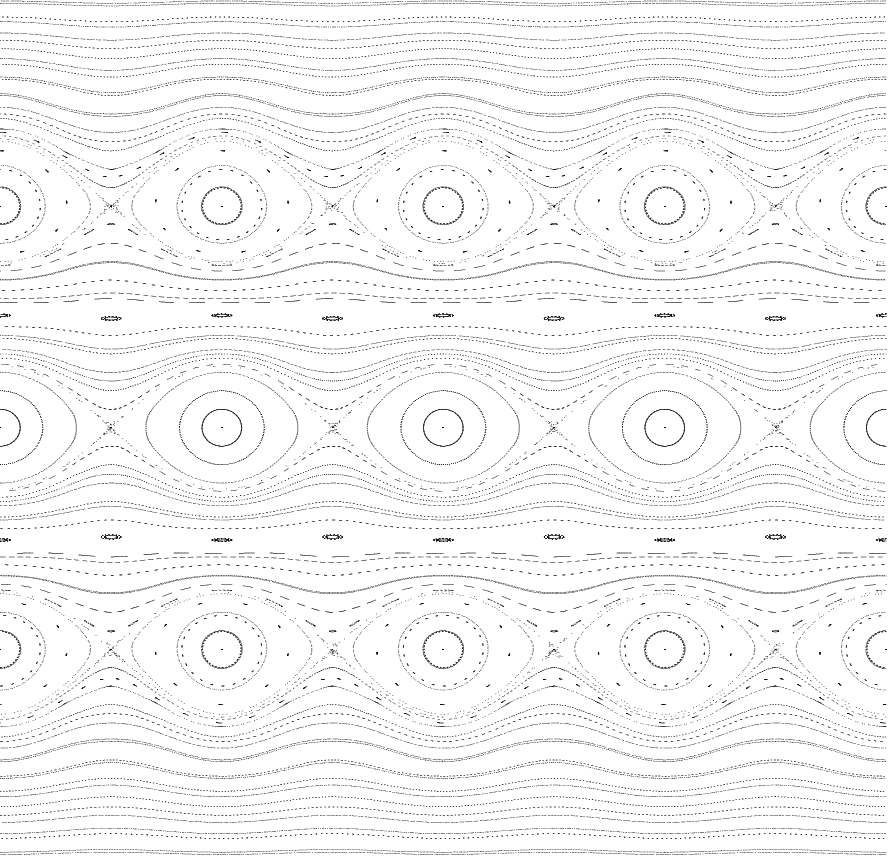} \hskip 1cm
		\includegraphics[height=4cm]{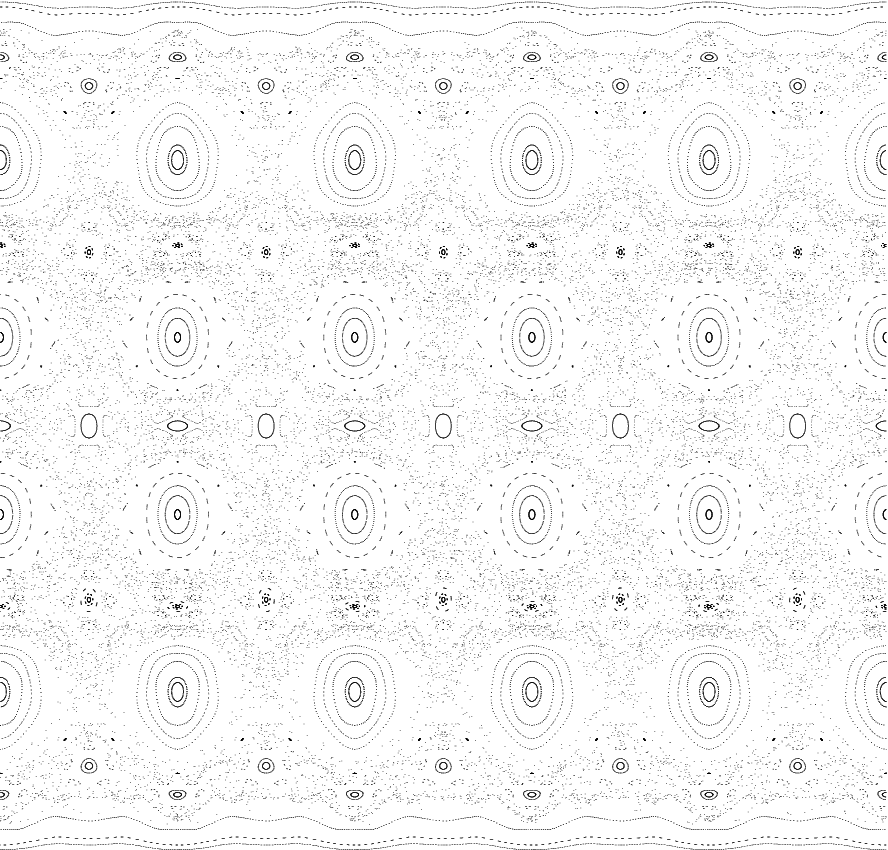}
	\end{center}    
	\caption{General structure of the phase space of symmetric billiards}
%(L) $4$-symmetric oval with $R(\phi)=1+ 0.2 \cos (4\phi) + 0.01 \cos (8\phi )$.
%(R) $5$-symmetric oval with $R(\phi)=1+ 0.4\cos (5\phi) + 0.25 \cos (10\phi)$.}
\end{figure}

\begin{figure}[h]
\label{fig:retas}
		\includegraphics[width=0.3\hsize]{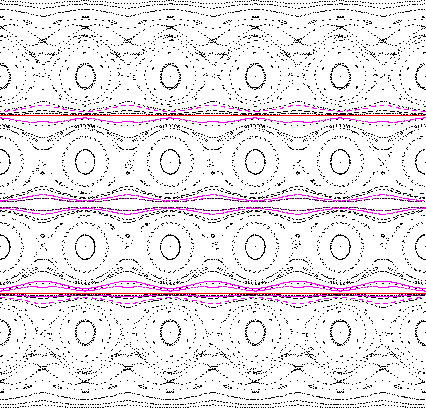} \hfill
		\includegraphics[width=0.3\hsize]{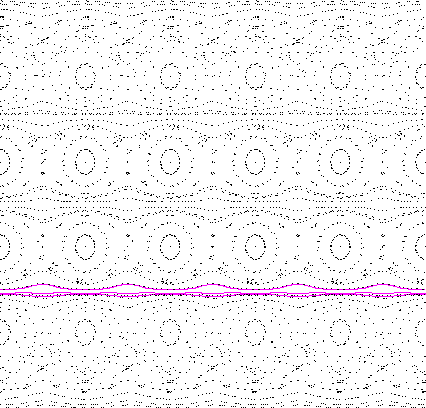} \hfill
		\includegraphics[width=0.3\hsize]{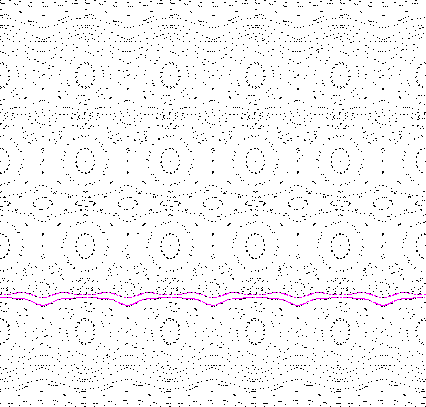} 
	\caption{Invariant curves}
\end{figure}

Tabachnikov \cite{tab} and Gutkin \cite{gutkin} proved that a non circular oval billiard has an horizontal invariant curve of the form $p=p_0 =\cos \alpha_0 \neq 0$ if and only if the radius of curvature of the boundary is given by $R(\varphi)=1+a_1\cos(n\varphi)$ (so the curve is symmetric), with $n\geq 4$, $|a_1|<1$, and $\alpha_0$ satisfies $\tan (n\alpha_0) = n\tan(\alpha_0)$. 
If $n$ is odd, $\alpha_0=\pi/2$ is a solution and the billiard's boundary has constant width. 
In this case, one can prove the existence of invariant curves close to $p_0=0$ \cite{tese}.
The persistence of invariant curves under perturbations of the boundary (as observed on Figure~\ref{fig:retas}), although expected, is a challenging subject.
%%%%%%%%%%%%%%%%%%%%%%%%%%%%%%%%%%%%%%%%%%%%%%%%%%%%%%%%%%%%

\begin{comment}
\begin{itemize}
\item Geraldo, voce provou que perturbacoes de curvas com largura constante vao ter curvas invariantes perto do p=0?
\item Citar o  resultado do Bunimovich (checar se  é genérico ou aberto e denso) ????
\item largura constante (simetria impar)
\item relacao entre regioes de instabilidade e as orbitas de simetricas
\item separacao de regioes de instabilidade
\item coisas conectadas onde aplicar: elliptic flowers, gutkin billiards
\end{itemize}
\end{comment}

Finally we observe that most of our results may be extended to piecewise smooth curves such as Bunimovich's flowers \cite{bunin22}.
%\notaIL{Syok}{curiosidade: como e h a funcao suporte de uma curva de largura constante?} %final

%%%%%%%%%%%%%%%%%%%%%%%%%%%%%%%%%%%%%%%%%%%%%%%%%%%%%%%
\begin{ack}
This work originated from the thesis of G.C.G. Ferreira \cite{tese} which
was financed by the Coordena\c c\~ao de Aperfei\c coamento de Pessoal de N\'\i vel Superior (CAPES) and 
Conselho Nacional de Desenvolvimento Cient\'\i fico e Tecnol\'ogico (CNPq),  Brasil.  We want to thank L. Bunimovich for his interest in this work.
\end{ack}
\bibliographystyle{plain}
\bibliography{sym-refs}

%\newpage
%\listoftodos\newpage
%\listoffigures
 \end{document}